\documentclass[12pt,leqno]{article}
\usepackage{amsmath,amsthm}

\def\init{\setcounter{equation}{0}}
\newtheorem{theorem}{Theorem}[section]

\newtheorem{prop}{Proposition}[section]
\newcommand{\R}{{\bf R}}
\newcommand{\C}{{\bf C}}
\newcommand{\Z}{{\bf Z}}

\newtheorem{lemma}{Lemma}[section]

\newcommand{\e}{{\varepsilon}}

\title{Inverse problems for the Schr\"{o}dinger  equations
with time-dependent electromagnetic potentials and
the Aharonov-Bohm effect.
\author{G.Eskin, \ \ \  Department of Mathematics, UCLA,\\ Los Angeles,
CA 90095-1555, USA. \ E-mail: eskin@math.ucla.edu}
}

\begin{document}

\maketitle
\begin{abstract}
We consider the inverse boundary value problem for the Schr\"{o}dinger operator
with time-dependent electromagnetic potentials in domains with obstacles.  We
extend the resuls of the author's works [E1], [E2], [E3]
to the case of time-dependent potentials.  We relate our results to the Aharonov-Bohm
effect caused by magnetic and electric fluxes.  
\end{abstract}

\section{Introduction.}
\label{section 1}
\init

Let $\Omega_0$ be a smooth  simply-conneted domain in $\R^n$  not necessarily bounded.
Let $D\subset\Omega_0\times [0,T]$  be a domain with the following 
properties : 
\\
Denote  $D_{t_0}=D\cap\{t=t_0\}$.  We assume that
$D_{t_0}=\Omega_0\setminus\cup_{j=1}^m\overline{\Omega_j(t_0)}$  where 
$\Omega_j(t_0)$ are a piece-wise smooth nonintersecting domains,
$1\leq j\leq m$.  We assume that $\Omega'(t_0)=\cup_{j=1}^m\Omega_j(t_0)$
depends smoothly on $t_0\in[0,T]$.  We do not assume that $\Omega_j(t_0)$
are bounded.

Finally we assume that the normal to $\partial D\setminus( D_0\cup D_T)$
in $\R^{n+1}$  is not parallel to the $t$-axis for any $t\in [0,T]$.

Consider the Schr\"{o}dinger equation with time-dependent electromagnetic potentials
in $D$:

\begin{equation}                          \label{eq:1.1}
i\frac{\partial u(x,t)}{\partial t} = Hu\stackrel{def}{=}\sum_{j=1}^n
(-i\frac{\partial}{\partial x_j}-A_j(x,t))^2u(x,t)+V(x,t)u(x,t),\ \ (x,t)\in D,
\end{equation}
with zero initial conditions
\begin{equation}                           \label{eq:1.2}
u(x,0)=0,\ \ \ \ (x,0)\in D_0,
\end{equation}
and the Dirichlet boundary conditions 
\begin{eqnarray}                          \label{eq:1.3}
u|_{\partial\Omega_0\times(0,T)}=f,
\\
u|_{\partial\Omega'(t_0)}=0,\ \ \forall t_0\in [0,T].
\nonumber
\end{eqnarray}

We assume that $A(x,t),V(x,t)$ are smooth in $D$ with  compact support.
Without loss of generality we can assume that $A(x,t)=0,V(x,t)=0$  near
$\partial\Omega_0\times[0,T]$  (see Remark 3.1).
Let $\Lambda$  be the Dirichlet-to-Neumann (D-to-N) operator on
$\partial\Omega_0\times[0,T]$, i.e.
\begin{equation}                         \label{eq:1.4}
\Lambda f 
=\frac{\partial u}{\partial \nu}-i(A\cdot\nu)u|_{\partial\Omega_0\times(0,T)},
\end{equation}
where $f$ is the same as in (\ref{eq:1.3}).  Denote by $G(D)$  the group of 
$C^\infty(\overline{D})$ functions $c(x,t)$  such that $c(x,t)\neq 0$
in $\overline{D}$ and denote by $G_0(D)$ the subgroup of $G(D)$ such that 
$c(x,t)=1$  on $\partial\Omega_0\times[0,T]$.  We say that 
the electromagnetic potentials $(A(x,t),V(x,t))$ and  $(A'(x,t),V'(x,t))$  
are gauge equivalent if there exists $c(x,t)\in G(D)$  such that
\begin{eqnarray}                               \label{eq:1.5}
A_j'(x,t)=A_j(x,t)+ic^{-1}(x,t)\frac{\partial c}{\partial x_j},\ \ \ 1\leq j\leq n,
\\
V'(x,t)=V(x,t)-ic^{-1}(x,t)\frac{\partial c}{\partial t}.
\nonumber
\end{eqnarray}
Note that if $(i\frac{\partial}{\partial t} - H)u=0$  and
\begin{equation}                           \label{eq:1.6}
u'=c^{-1}u
\end{equation}
then $(i\frac{\partial}{\partial t}-H')u'=0$,  where $H'$ is the operator
with potentials $A'(x,t),V'(x,t)$.  The 
group $G(D)$ is called the gauge group and (\ref{eq:1.6})  is called
the gauge transformation.   When  $A(x,t)$  and $V(x,t)$
are real-valued it is natural to consider only $c(x,t)$ such that
$|c(x,t)|=1$.
Also when $D$ is simply connected any $c(x,t)\in G(D)$ has a form 
$c(x,t)=e^{i\varphi(x,t)}$,  where  $\varphi(x,t)\in C^\infty(\overline{D})$.

We say that the D-to-N operators $\Lambda,\ \Lambda'$ corresponding to
the Schr\"{o}dinger operators $i\frac{\partial}{\partial t}-H,\ 
i\frac{\partial}{\partial t}-H'$  are gauge equivalent if there exists 
$c(x,t)\in G(D)$  such that 
\begin{equation}                      \label{eq:1.7}
\Lambda'=c_0^{-1}\Lambda c_0\ \ \ \mbox{on}\ \ \partial\Omega_0\times[0,T],
\end{equation}
where $c_0$  is the restriction of $c(x,t)$  to $\partial \Omega_0\times[0,T]$.

We shall introduce gauge invariant boundary data on 
$\partial\Omega_0\times(0,T)$  (c.f.  [E1]).
Let  $u(x,t)$  be a solution of $i\frac{\partial u}{\partial t}- H u =0$
in $D,\ u(x,0)=0.$

Let $|u(x,t)|^2$ be te probability density  and 
\[
S(x,t)=\Im\left(\frac{\partial u(x,t)}{\partial x}
-iA(x,t)u(x,t)\right)\overline{u(x,t)}
\] 
be the probability current.
We define the boundary data of $u(x,t)$ on $\partial\Omega_0\times (0,T)$
as 
\begin{equation}                              \label{eq:1.8}
|u(x,t)|^2|_{\partial\Omega_0\times(0,T)}=f_1,\ \
\frac{\partial}{\partial \nu}|u(x,t)|^2|_{\partial\Omega_0\times(0,T)}=f_2,\ \ 
S(x,t)|_{\partial\Omega_0\times(0,T)}=f_3,
\end{equation}
where
$\frac{\partial}{\partial \nu}$ is the normal derivative.  Note  that 
the probability density and the probability current are gauge invariant.

The following proposition was proven in [E1]:

\begin{prop}                         \label{prop:1.1}                          
  Let
$(i\frac{\partial}{\partial t}-H)u=0$  and 
$(i\frac{\partial}{\partial t}-H')u'=0$
be two Schr\"{o}dinger equations in $D,\ u(x,0)=u'(x,0)=0,
\ u(x,t)|_{\partial\Omega'(t)}=
u'|_{\partial\Omega'(t)}=0$   for all $t\in [0,T]$.   Then the D-to-N
operators $\Lambda$  and $\Lambda'$  are gauge equivalent on 
$\partial\Omega_0\times(0,T)$  if and only if for any solution $u$ there exists 
$u'$ such that the gauge invariant boundary data (\ref{eq:1.8}) of
$u(x,t)$ and $u'(x,t)$ are equal.
\end{prop}

This paper is 
a completion  of works [E1], [E2], [E3]  on the inverse problems for 
the stationary Schr\"{o}dinger  equation in a situation where the Aharonov-Bohm
effect caused by magnetic fluxes holds.  Here we extend these results 
to the time-dependent Schr\"{o}dinger equation with time-dependent electromagnetic
potentials in the cases when the Aharonov-Bohm effect holds caused by both the
magnetic and electric fluxes.

Consider two Schr\"{o}dinger operators $i\frac{\partial}{\partial t}-H^{(p)},\ p=1,2,$
with time-dependent electromagnetic potentials  $(A^{(p)}(x,t),V^{(p)}(x,t)),\ p=1,2,$ 
in the domain $D$ described above.

\begin{theorem}                                      \label{theo:1.1}
Assume that $n\geq 2$  and the domain $D$ is such that for any 
$t_0\in [0,T]\ \ \ D_{t_0}$ satisfies the following condition (c.f.  [E2],  page 287):
All $\Omega_j(t_0)$ are convex,  $1\leq j\leq m$,  and for each point 
$x_0\in D_{t_0}$  there exists a two-dimensional plane $\Pi\subset \R^n$
such that $\Pi$  intersects at most one domain $\Omega_j(t_0),\ 1\leq j\leq m$.
 Let $i\frac{\partial}{\partial t}-H_j,j=1,2,$  be two Schr\"{o}dinger operators
in $D$ and $\Lambda^{(j)},j=1,2,$  be the corresponding D-to-N operator on
$\partial\Omega_0\times(0,T)$.  If  $\Lambda^{(1)}$  and $\Lambda^{(2)}$ 
are gauge equivalent on $\partial\Omega_0\times(0,T)$  then the electromagnetic
potentials $A^{(1)}(x,t), V^{(1)}(x,t)$  and $A^{(2)}(x,t),V^{(2)}(x,t)$
are gauge equialent in $D$.
\end{theorem}

Theorem \ref{theo:1.1} generalized the result of [E2].
Note that when $n=2$  the domain  $D_{t_0}$  contain at most  one
convex $\Omega_j(t_0)$.

The  following theorem generalizing the result of [E3]  treats  the case 
when $D_{t_0}$ contains more than one obstacle $\Omega_j(t_0)$:
\begin{theorem}                            \label{theo:1.2}
Assume that $n=2$  and that for each $t_0\in[0,T]$ domains 
$\Omega_j(t_0),\ 1\leq j\leq m$,
are  piece-wise smooth and convex.  Suppose also that 
there is no trapping broken rays (more exactly,  the conditions 
a),  b) in [E3],  page 1507,  are satisfied).    If D-to-N
operators $\Lambda^{(1)}$ and $\Lambda^{(2)}$  are gauge equivalent on
$\partial\Omega_0\times(0,T_0)$, then the electromagnetic potentials 
$A^{(1)}(x,t),V^{(1)}(x,t)$ and $A^{(2)}(x,t),V^{(2)}(x,t)$
are gauge equivalent on $D$.
\end{theorem}

Inverse problems and Aharonov-Bohm effect were considered in [N],  [W].
General inverse problems for the stationary Schr\"{o}dinger equations can
be reduced to the inverse problems for the hyperbolic equations with the 
time-independent coefficients (see [B], [KKL], [E6]  and additional 
references there). 
Hyperbolic  equations approach is very powerful.   However it does not 
apply to the Schre\"{o}dinger equations with time-dependent potentials. 

The plan of the paper is the following:
In \S 2 we prove  Theorems \ref{theo:1.1} and \ref{theo:1.2}.
In \S 3 we consider the Schr\"{o}dinger operators with time-dependent
Yang-Mills potentials.
In \S 4 we discuss the Aharonov-Bohm effect (see [AB], [WY], [OP])  in
the situations covered by Theorems \ref{theo:1.1}   and \ref{theo:1.2}.
Note that the Aharonov-Bohm effect is caused by the magnetic and
electric fluxes,  and the cause of electric fluxes appears only  when
potentials are time-dependent.

\section{The proof of Theorems \ref{theo:1.1}  and \ref{theo:1.2}.}
\label{section 2}
\init

We shall start with the construction of geometric optics solutions of (\ref{eq:1.1}).
We are looking for a solution in the form
\begin{equation}                            \label{eq:2.1}
u_{N0}=e^{-ik^2t+ik(x\cdot \omega)}a_0^{(N)},
\end{equation}
where $a_0^{(N)}=\sum_{p=0}^N\frac{a_{p0}(x,t,\omega)}{(ik)^p},\ |\omega|=1,\ k$
is a large parameter.  Substituting
$u_{N0}$ in (\ref{eq:1.1}) we obtain 
\[
k^2a_0^{(N)}+i\frac{\partial a_0^{(N)}}{\partial t} =
(-i\frac{\partial}{\partial x}+k\omega- A(x,t))^2a_0^{(N)}+V(x,t)a_0^{(N)}.
\]
Equating the equal power of $k$ we get
\begin{equation}                            \label{eq:2.2}
\omega\cdot(-i\frac{\partial}{\partial x} -A)a_{00}=0,
\end{equation}
\begin{equation}                            \label{eq:2.3}
2\omega\cdot(-i\frac{\partial}{\partial x} -A)a_{p0}=
(i\frac{\partial}{\partial t} -H)a_{p-1,0},\ \ p\geq 1.
\end{equation}
To solve (\ref{eq:2.2}) denote $s=\omega\cdot x, \ \tau_j=(\omega_{\perp j},x),
\ 1\leq j\leq n-1$,  where $\omega_{\perp j},1\leq j\leq n-1,$   
is a basis in an orthogonal complement to $\omega$ in $\R^n$.
We have $\frac{\partial a_{00}}{\partial s}-i(A\cdot \omega)a_{00}(s,\tau,t)=0$,
where $\tau=(\tau_1,...,\tau_{n-1})$.  Therefore we can take 
\begin{equation}                            \label{eq:2.4}
a_{00}(s,\tau,t)=\chi_1(t)\chi_2(\tau)
\exp\left(i  \int_{s_0}^s
A(\tau+s'\omega)\cdot\omega ds'\right) ,
\end{equation}
where 
$\chi_1(t)=\frac{1}{\sqrt{\e}}\chi_0(\frac{t-t_0}{\e}),\ 
\chi_2(\tau)
=\frac{1}{\e^{\frac{n-1}{2}}}\Pi_{j=1}^{n-1}\chi_0(\frac{\tau_j-\tau_{0j}}{\e}),\ 
\chi_0\in C_0^\infty(\R^1),\ \chi_0(t)=0$ for 
$|t|>1,\ \int_{-\infty}^\infty\chi_0^2(t)dt=1$.

We solve (\ref{eq:2.3}) prescribing the initial condition
\begin{equation}                         \label{eq:2.5}
a_{p0}(s_0,\tau,t)=0,\ \ 1\leq p\leq N.
\end{equation}

Here $s_0,\tau_0,t_0$ are such that
the point $(x_0^{(0)}=s_0\omega+\sum_{j=1}^{n-1}\tau_{0j}\omega_{\perp j},t_0)$
does not belong
to $\Omega_0\times[0,T]$.  We have 
\[
(i\frac{\partial}{\partial t}-H)u_{N0}=f_{N0},
\]
where $\hat{f}_{N0}=\frac{e^{-ik^2t+ik(\omega\cdot x)}}{(ik)^N}
(i\frac{\partial}{\partial t}-H)a_{N0}.$

Until now we assume that the ray $x=(s-s_0)\omega+x_0^{(0)}$ 
 does not meet the obstacles. 
In the case when this ray meets an obstacle we have to consider a broken ray
that reflects at $\partial\Omega_j(t_0),\ 1\leq j\leq m$ (c.f.  [E3],  pages
1498-99).  More precisely,  we say that $\gamma=
\gamma_0\cup\gamma_1\cup...\cup\gamma_r$
is a broken ray with legs $\gamma_0,\gamma_1,...,\gamma_r$  if $\gamma_0$  is 
the ray $x=x_0^{(0)}+(s-s_0)\omega,\ s_0\leq s\leq s_1,\
\gamma_0$  hits
$\partial\Omega'(t_0)=\cup_{j=1}^m\partial\Omega_j(t_0)$ at some point
$x_0^{(1)}=x_0^{(0)}+(s_1-s_0)\omega\in \partial\Omega'(t_0),\ \gamma_1$
is the reflected ray.  We assume that
$x_0^{(1)}$ is not a tangential point of reflection.  Then
\begin{equation}                           \label{eq:2.6}
\theta_1=\theta_0-2(n_1(x_0^{(1)},t_0)\cdot\theta_0)n_1(x_0^{(1)},t_0),
\end{equation}   
where
$\theta_0=\omega$  is
the direction of $\gamma_0,\theta_1$  is the direction of $\gamma_1,n(x_0^{(1)},t_0)$
is the outward unit normal to $\partial\Omega'(t_0)$  at the point $x_0^{(1)}$.  
Analoguosly,  $\gamma_2,...,\gamma_{r-1}$ have  nontangential points of reflection
on $\partial\Omega'(t_0)$  and $\gamma_r$  ends on $\partial\Omega_0\times\{t=t_0\}$.

We associate with broken ray $\gamma$  the following geometric optics solution 
(c.f. [E3],  page 1499):
\begin{equation}                    \label{eq:2.7}
u_N=\sum_{j=0}^r\sum_{p=0}^N\frac{a_{pj}(x,t,\omega)}{(ik)^p}
e^{-ik^2t+ik\psi_j(x,t,\omega)},
\end{equation}
where  $\psi_0(x,t,\omega)=x\cdot\omega,\ a_{0p}(x,t,\omega)$  are the same as in 
(\ref{eq:2.2}), (\ref{eq:2.3}), (\ref{eq:2.4}), (\ref{eq:2.5}),
eiconals  $\psi_j(x,t,\omega)$   satisfy  the equations
\begin{eqnarray}                       \label{eq:2.8}
\left|\frac{\partial \psi_j(x,t,\omega)}{\partial x}\right|=1,\ \ \ \ \ \ \ \ \ \ 
\ \ \ \ \ \ \ \ \ 
\\
\psi_j(x,t,\omega)|_{\partial\Omega'(t)}=\psi_{j+1}(x,t,\omega)|_{\partial\Omega'(t)},
\\
\frac{\partial\psi_{j+1}(x_0^{(j+1)},t_0,\omega)}{\partial x}=\theta_{j+1},\ \ 
0\leq j\leq r-1,
\end{eqnarray}
where
$x_0^{(j+1)}$  is the point of reflection of $\gamma_j$ at $\partial\Omega'(t_0)$ and
$\theta_{j+1}$ is the direction of $\gamma_{j+1}$.
We assume that $\gamma_j,1\leq j\leq r$,  do not contain caustic points.
Then the eiconals $\psi_j,1\leq j\leq r,$ exist.  In particular,  this is true
when  $\Omega_j(t)$ are convex.  Functions $a_{pj}$ satisfy the following equations:
\begin{equation}                                       \label{eq:2.11}
2\frac{\partial a_{pj}}{\partial x}\cdot\frac{\partial\psi_j}{\partial x}
+\Delta\psi_j a_{pj} -2iA(x,t)\cdot\frac{\partial\psi_j}{\partial x}a_{pj}=
f_{pj}(x,t,\omega)+i\frac{\partial \psi_j}{\partial t}a_{pj},\ p\geq 0,
\end{equation}
where  $f_{0j}=0,\ f_{pj}$ depends on $a_{0j},...,a_{p-1,j}$.
When $\Omega'(t)$  is independent of $t$  then $\psi_j,\ j\geq 1$  is also
independent of $t$ but  $\psi_j(x,t,\omega)$  depends on $t$  when
$\Omega'(t)$   depends on $t$.
We impose the following conditions on $a_{pj}$:
\begin{equation}                      \label{eq:2.12}
a_{pj}|_{\partial\Omega'(t)}=-a_{p,j+1}|_{\partial\Omega'(t)},
\ \ \ 0\leq j\leq r-1.
\end{equation}
This conditions imply that 
\begin{equation}                     \label{eq:2.13}
u_N|_{\partial\Omega'(t)}=0.
\end{equation}
Substituting  (\ref{eq:2.7})
into (\ref{eq:1.1})  we get
\begin{equation}                    \label{eq:2.14}
(i\frac{\partial}{\partial t}-H)u_N=f_N,
\end{equation}
where  $f_N=\sum_{j=0}^r\frac{f_{Nj}}{(ik)^N}e^{ik\psi_j-ik^2t}\ \ f_{Nj}(x,t,\omega)$
are smooth.

In order to 
complete the construction of the geometric optics solutions
we need the following lemma:
\begin{lemma}                          \label{lma:2.1}
Consider the initial-value problem
\begin{eqnarray}                       \label{eq:2.15}
(i\frac{\partial}{\partial t}-H)w=f \ \ \mbox{in\ \ } D,
\\
w|_{\partial D_t}=0, \ \ \ \forall t\in (0,T),
\nonumber
\\
w(x,0)=w_0(x),\ \ (x,0)\in D_0.
\nonumber
\end{eqnarray}

Suppose that $\Omega'(t)$  is a smooth  domain  in $\R^n$ for
$\forall t\in [0,T]$.
Then for
any $f\in C^1([0,T],L^2(D_t))$  and any 
$w_0\in H^2(D_0)\cap\stackrel{\circ}{H^1}(D_0)$
there exists a unique $w\in C([0,T],\stackrel{\circ}{H^1}(D_t)\cap H^2(D_t))\cap C^1([0,T],L^2(D_t))$.
\end{lemma}

{\bf Proof:}
We assume that $A(x,t)$ are real but $V(x,t)$ can be complex-valued.

Make change of variables:
\begin{equation}                             \label{eq:2.16}
t=t,\ x'=\varphi(x,t),   \ \ 
\varphi(x,t)=(\varphi_1(x,t),...,\varphi_n(x,t))
\end{equation}
in $\overline{D}$ such that $\varphi(x,t)=x$  on $\partial\Omega_0\times[0,T]$  and
$\varphi(\Omega_j(t),t)=\Omega_j(0),\ \forall t\in[0,T]$.
Denote $\hat{u}(x',t)=u(x,t)$,  where $x'=\varphi(x,t)$.  Then 
$u_t(x,t)=\hat{u}_t+\hat{u}_{x'}(x',t)\cdot\varphi_t(x,t),\ 
u_x(x,t)=u_{x'}\frac{\partial \varphi}{\partial x}$, etc.

Therefore (\ref{eq:2.15}) has the following form in $(x',t)$ coordinates:
\begin{eqnarray}                           \label{eq:2.17}
i\frac{\partial \hat{u}}{\partial t}=
\hat{H}_0\hat{u}+\hat{H}_1\hat{u} +\hat{f},\ \ (x',t)\in \hat{D}
\\
\hat{u}|_{\partial\hat{\Omega}\times(0,T)}=0,\ \ \hat{u}(x',0)=
\hat{u}_0(x'),\ x'\in\hat{\Omega},
\nonumber
\end{eqnarray}
where 
\begin{eqnarray}
\hat{H}_0=\frac{1}{\sqrt{g(x',t)}}\sum_{j,k=1}^n(-i\frac{\partial}{\partial x_j'}
-\hat{A}_j(x',t))\sqrt{g}g^{jk}(x',t)(-i\frac{\partial}{\partial x_k'}-\hat{A}_k),
\nonumber
\\
H_1=\hat{V}(x',t)-i\varphi_t\cdot\hat{u}_{x'},\ \ \ \ \ \ \ \ \ \ \ \ \ \ \ \ \ 
\ \ \ \ \ \ \ \ \ \ \ \ \ \ \ 
\nonumber
\\
g^{jk}(x',t)=\sum_{p=1}^n
\frac{\partial \varphi_j}{\partial x_p}\frac{\partial \varphi_k}{\partial x_p},\ \ g(x',t)=\det\|g^{jk}\|^{-1}.\ \ \ \ \ \ \ \
\nonumber
\end{eqnarray}
Note $\hat{D}=\hat{\Omega}\times(0,T],\ \ \hat{\Omega}
=\Omega_0\setminus\overline{\Omega}'$,  where $\Omega'=
\cup_{j=1}^m\Omega_j(0)$.

Multiplying (\ref{eq:2.17}) by $\sqrt{g(x',t)}\ \hat{u}(x',t)$
and integrating over $\hat{\Omega}$  we get 
\begin{equation}                                \label{eq:2.18}
(i\frac{\partial\hat{u}}{\partial t},\hat{u})_g=
(\hat{H}_0\hat{u}.\hat{u})_g+(\hat{H}_1u+f,\hat{u})_g,
\end{equation}
where $(v_1,v_2)_g$ is the $L_2$ inner product in $\hat{\Omega}$ with the weight
$\sqrt{\hat{g}}$.  Take the imaginary part of (\ref{eq:2.18}) and
integrate in $t$ from $0$ to $t$.   Since $\hat{H}_0$ is self-adjoint and since 
\begin{equation}                             \label{eq:2.19}
\int_{\hat{\Omega}}i\sum_j\varphi_{jt}\frac{\partial\hat{u}}{\partial x_j}
\sqrt{g}\overline{\hat{u}}dx'=
-\int_{\hat{\Omega}}i\sum_j\varphi_{jt}\sqrt{g}\hat{u}
\overline{\frac{\partial\hat{u}}{\partial x_j}}dx'+
\int_{\hat{\Omega}} b(x',t)|\hat{u}|^2dx'
\end{equation}
for some $b(x',t)$,  we get
\begin{eqnarray}                              \label{eq:2.20}
\int_{\hat{\Omega}}\sqrt{g(x',t)}|\hat{u}(x',t)|^2dx'-\int_{\hat{\Omega}}
\sqrt{g(x',0)}|\hat{u}_0(x',0)|^2dx'
\\
\leq C\int_0^t\|\hat{u}\|_0^2dt'+
C\int_0^t\|\hat{f}\|_0\|\hat{u}\|_0dt',
\nonumber
\end{eqnarray}
where 
$\|\ \ \|_0$ is the $L_2$ norm in $\hat{\Omega}$.

It follows from (\ref{eq:2.20})
that
\begin{equation}                           \label{eq:2.21}
\max_{[0,T]}\|\hat{u}(\cdot,t)\|_0^2\leq C\|\hat{u}(x',0)\|_0^2+
C\int_0^T\|\hat{f}(\cdot,t)\|_0^2dt'.
\end{equation}
  Changing from the beginning in (\ref{eq:1.1})
$u$ to $e^{i\lambda t}u$  we can assume that  $\Re V(x,t)$ is  large.
Since $\hat{u}|_{\partial\hat{\Omega}\times (0,T)}=0$  and 
$\hat{H}_0+\hat{H}_1$  is elliptic for each $t\in[0,T]$ we get from (\ref{eq:2.17}),
 by the standard
elliptic theory:
\begin{equation}                          \label{eq:2.22}
\|\hat{u}\|_2\leq C\left\|\frac{\partial\hat{u}}{\partial t}\right\|_0
+C\|\hat{f}\|_0,
\end{equation}
where
$\|\hat{u}\|_2$  is the Sobolev norm in $\hat{\Omega}$.

Differentiate (\ref{eq:2.17}) in $t$:
\begin{equation}                         \label{eq:2.23}
i\frac{\partial^2\hat{u}}{\partial t^2}
=\hat{H}_0\frac{\partial\hat{u}}{\partial t}
+\hat{H}_1\frac{\partial\hat{u}}{\partial t}+
\hat{H}'\hat{u}+\frac{\partial \hat{f}}{\partial t},
\end{equation}
where $\hat{H}'$ is a differenial operator in $x$ of order at most 2.  
Note that
$\frac{\partial\hat{u}}{\partial t}|_{\partial\hat{\Omega}\times [0,T]}=0$
since $\hat{u}|_{\partial\hat{\Omega}\times [0,T]}=0$.
Multiplying (\ref{eq:2.23})  by $\sqrt{g}\overline{\frac{\partial\hat{u}}
{\partial t}}$,
integrating over $\hat{\Omega}\times (0,t)$  and taking  the imaginary
part we get,  as in (\ref{eq:2.20}):
\begin{eqnarray}                           \label{eq:2.24}
\int_{\hat{\Omega}}\sqrt{g}\left|\frac{\partial\hat{u}}{\partial t}\right|^2dx'
-\int_{\hat{\Omega}}\sqrt{g(x',0)}\left|\frac{\partial\hat{u}(x',0)}{\partial t}\right|^2dx'
\\
\leq C\int_0^t\left\|\frac{\partial\hat{u}}{\partial t}\right\|_0^2dt'
+ C\int_0^t\|\hat{u}\|_2\left\|\frac{\partial\hat{u}}{\partial t}\right\|_0dt'
+ C\int_0^t\left\|\frac{\partial\hat{f}}{\partial t}\right\|_0
\left\|\frac{\partial\hat{u}}{\partial t}\right\|_0dt'.
\nonumber
\end{eqnarray}
Using
(\ref{eq:2.22})  we get from  (\ref{eq:2.24})
\begin{equation}                           \label{eq:2.25}
\max_{[0,T]}\left\|\frac{\partial\hat{u}}{\partial t}\right\|_0^2
\leq C \left\|\frac{\partial\hat{u}(\cdot,0)}{\partial t}\right\|_0^2
+C\int_0^T\left(\|\hat{f}\|_0^2+
\left\|\frac{\partial\hat{f}}{\partial t}\right\|_0^2\right)dt'.
\end{equation}
It follows from (\ref{eq:2.17}) that 
\begin{equation}                         \label{eq:2.26}
\left\|\frac{\partial\hat{u}(x',0)}{\partial t}\right\|_0
\leq C\|\hat{u}_0(x',0)\|_2+\|\hat{f}(x',0)\|_0
\end{equation}
Combining (\ref{eq:2.21}), (\ref{eq:2.22}),  (\ref{eq:2.25}),  (\ref{eq:2.26})
we get 
\begin{equation}                          \label{eq:2.27}
\max_{[0,T]}\|\hat{u}\|_2^2+\max_{[0,T]}\left\|\frac{\partial\hat{u}}{\partial t}
\right\|_0^2\leq C\|\hat{u}_0\|_2^2+C\max_{[0,T]}\|\hat{f}\|_0^2+C\int_0^T\left\|
\frac{\partial\hat{f}}{\partial t}\right\|_0^2dt'
\end{equation}
We proved estimate (\ref{eq:2.27})  assuming  that the solution 
$\hat{u}$ exists.  To prove the existence 
we shall use the parabolic regularization.

Consider the parabolic equation:
\begin{equation}                          \label{eq:2.28}
(i-\e)\frac{\partial u_\e}{\partial t}=
\hat{H}_0u\e +\hat{H}_1u_\e +\hat{f},\ \ \e>0,
\end{equation}
with the same boundary and initial conditions as in (\ref{eq:2.17}):
\begin{equation}                          \label{eq:2.29}
u_\e|_{\partial\hat{\Omega}\times(0,T)}=0,\ \ u_\e(x',0)=\hat{u}_0(x'),
\end{equation}
By the parabolic theory there exists a solution of (\ref{eq:2.28}),
(\ref{eq:2.29})
belonging to the same space as $\hat{u}$.  Note that the estimate 
(\ref{eq:2.27})
holds for $u_\e$ with constants in.dependent  of  $\e>0$.  Taking the
weak limit of a subsequence $\e_k\rightarrow 0$  we get the existence
of $\hat{w}$ satisfying (\ref{eq:2.17}).   Note that in the original coordinates
$x=\varphi^{-1}(x',t)$ the estimate (\ref{eq:2.27})  still holds.
\qed

{\bf Remark 2.1}.  In the case of several obstacles we require that
$\Omega_j(t),\ 1\leq j\leq m$,  are convex sets having corners.  In this case
the estimates (\ref{eq:2.22})  does not hold and one has a weaker estimate:
\begin{equation}                              \label{eq:2.30}
\|\hat{u}\|_1\leq C\left\|\frac{\partial \hat{u}}{\partial t}\right\|_0
+C\|\hat{f}\|_0.
\end{equation}
To prove the existence  of the solution of (\ref{eq:2.17}) in this case we assume
for simplicity  that obstacles $\Omega_j(t)$  move as rigid  bodies,  i.e.
the change of variables (\ref{eq:2.16})  has the form 
$t=t,\ \varphi(x,t)=x+\varphi_j(t),\ \varphi_j(0)=0$, near  
$\Omega_j(t),\ 1\leq j\leq m$.  Then $g^{jk}(x',t)=\delta_{jk}$  near
$\Omega'$  and the differential operator $\hat{H}'$  has the  order
1  near $\Omega'$.  Therefore (\ref{eq:2.24})  holds with 
$\|\hat{u}\|_2$  replaced by  $\|\hat{u}\|_1+\|\chi u\|_2$  where $\chi=0$
near $\Omega'$.  Using the elliptic estimate of the form
(\ref{eq:2.22})  for $\chi\hat{u}$  and  the estimate (\ref{eq:2.30})
we get  that (\ref{eq:2.25})
holds in the case of  corners too.

Using the parabolic regulation as in (\ref{eq:2.28})
we get that there exists a unique solution  of (\ref{eq:2.15})
satisfying  the estimate (\ref{eq:2.27})  with $\|\hat{u}\|_2$  
replaced by $\|\hat{u}\|_1$.

By Lemma \ref{lma:2.1}  and Remark 2.1  there exists  $u^{(N+1)}\in 
C([0,T],H^1(D_t))$  such that
\begin{eqnarray}                          \label{eq:2.31}
\left(i\frac{\partial}{\partial t}-H\right)u^{(N+1)}=-f_N \ \ \ \mbox{in\ \ } D,
\\
u^{(N+1)}|_{\partial\Omega_j'(t)}=0,\ \ \ 1\leq j\leq m,\ \ \ 
\forall t\in [0,T],
\nonumber
\\
u^{(N+1)}(x,0)=0,\ \ \ u^{(N+1)}|_{\partial\Omega_0\times(0,T)}=0,
\nonumber
\end{eqnarray}
where  $f_N$  is the same as in (\ref{eq:2.14}).
Then  
$\|u^{(N+1)}\|_1+\left\|\frac{\partial u^{(N+1)}}{\partial t}\right\|_0
\leq Ck^{-N+1}$  for  $\forall\ t\in[0,T]$.
Therefore $u=u_N+u^{(N+1)}$  solves 
\begin{eqnarray}                                \label{eq:2.32}
\left(i\frac{\partial}{\partial t}-H\right)u=0,\ \ \ (x,t)\in D,
\ \ \
u|_{\partial\Omega_0\times(0,T)}=u_N|_{\partial\Omega_0\times(0,T)},
\nonumber
\\
u|_{\partial\Omega'(t)}=0,\ \  \forall t\in [0,T],\ \ u(x',0)=0,\ \ 
(x',0)\in D_0.
\end{eqnarray}

Now we shall use the Green's
formula.
Suppose we are given two time-dependent Schr\"{o}dinger operators 
$i\frac{\partial}{\partial t}-H_1$  and $i\frac{\partial}{\partial t}-H_2$  
such that the corresponding D-to-N operators 
$\Lambda_1$  and $\Lambda_2$ are gauge equivalent,    i.e. 
$\Lambda_1=g_0^{-1}\Lambda_2g_0$   where $g_0=g|_{\partial\Omega_0\times [0,T]},
g\in G(D)$.  Denote $A^{(3)} = A^{(2)}-ig^{-1}\frac{\partial g}{\partial x},\ 
V^{(3)}=V^{(2)}+ig^{-1}\frac{\partial g}{\partial t}$.  Let 
$i\frac{\partial}{\partial t}-H_3$  be the Schr\"{o}dinger
operator with electromagnetic potentials $A^{(3)},V^{(3)}$.  Then 
$i\frac{\partial}{\partial t}-H_1$  and $i\frac{\partial}{\partial t}-H_3$
have the same D-to-N operator : $\Lambda_3=\Lambda_1$.

Let $u_1$  be the solution of the initial-boundary value problem:
\begin{eqnarray}                                 \label{eq:2.33}
i\frac{\partial u_1}{\partial t}-H_1u_1=0 \ \ \ \mbox{in\ \ } D,
\\
u_1(x,0)=0,\ \ (x,0)\in D_0,
\nonumber
\\
u_1|_{\partial\Omega_0\times(0,T)}=f,\ \ \ u_1|_{\partial\Omega'(t)}=0
\ \ \ \forall t\in [0,T],
\nonumber
\end{eqnarray}
and let $u_3$  be the solution of     
\begin{eqnarray}                                 \label{eq:2.34}
i\frac{\partial u_3}{\partial t}-H_3^*u_3=0 \ \ \ \mbox{in\ \ } D,
\\
u_3(x,T)=0,\ \ (x,T)\in D_T,
\nonumber
\\
u_3|_{\partial\Omega_0\times(0,T)}=g,\ \ \ u_{3}|_{\partial\Omega'(t)}=0
\ \ \ \forall t\in [0,T],
\nonumber
\end{eqnarray}
where $H_3^*$  is the adjoint to $H_3$.

We have,  by the Green's formula (c.f. [SU]):
\begin{eqnarray}                             \label{eq:2.35}
0=((i\frac{\partial u}{\partial t}-H_1)u_1,u_3)-
(u_1,(i\frac{\partial}{\partial t}-H_3^*)u_3)
\\
=
\int_D(2iA^{(1)}(x,t)\cdot\frac{\partial u_1}{\partial x}\overline{u_3}
+2iu_1A^{(3)}\cdot\frac{\partial \overline{u_3}}{\partial x}
\nonumber
\\
+(q_1(x,t)-q_3(x,t))u_1\overline{u_3})dxdt+[\Lambda_1f,g]-
[f,\Lambda_3^*g],
\nonumber
\end{eqnarray}
where 
$q_p=i\frac{\partial}{\partial x}\cdot A^{(p)}(x,t)+A^{(p)}(x,t)\cdot
A^{(p)}(x,t)+V^{(p)}(x,t),\ [\ ,\ ]$  is the inner
product on $\partial\Omega_0\times[0,T]$ and $\Lambda_p$  are the D-to-N
operators corresponding to $H_p,\ p=1$ or 3.
Let
\begin{equation}                              \label{eq:2.36}
v=\sum_{j=0}^r\sum_{p=0}^N\frac{b_{pj}}{(ik)^p}e^{-ik^2t+ik\psi_j(x,t,\omega)}
+v^{(N+1)}
\end{equation}
be the geometric optics solution of (\ref{eq:2.34})
corresponding to the same broken ray $\gamma$  as the solution 
(\ref{eq:2.32}).  Substitute (\ref{eq:2.32}) in (\ref{eq:2.35})
instead of $u_1$  and substitute (\ref{eq:2.36})
instead of $u_3$.
Dividing (\ref{eq:2.35})  by $2k$,   passing to the limit when
$k\rightarrow\infty$  and taking
into account that $\Lambda_1=\Lambda_3$  we get (c.f. [E3], pages
1502-03):
\[
\sum_{j=0}^r\int_D(A^{(3)}(x,t)-A^{(1)}(x,t))\cdot\psi_{jx}(x,t,\theta)a_{0j}(x,t,\omega)
\overline{b_{0j}(x,t,\omega)})dxdt
\]
Note that $a_{00}$  and $b_{00}$  have the form (\ref{eq:2.4}).
Note also that terms containing $\psi_{jt}$  in $a_{0j}\overline{b_{0j}}$
cancel each other.
Making changes of variables  as in [E3], pages  1503-06  (see also [E4],  page 30),
and taking the limit when $\e\rightarrow 0$  where $\e>0$  is the same as in
(\ref{eq:2.4}),  we obtain
\begin{equation}                               \label{eq:2.37}
\exp[i\sum_{j=0}^r\int_{\gamma_j}(A^{(3)}(x_0^{(j)}+s\theta_j,t_0)
-A^{(1)}(x_0^{(j)}+s\theta_j,t_0))\cdot \theta_jds]=1,
\end{equation}
where $\theta_j$  is the direction of $\gamma_j,\ x_0^{(j)}$ is the starting  point of 
$\gamma_j,\ 0\leq j\leq r$.

It is important to emphasize that the ray $\gamma$  is contained in the plane 
$t=\mbox{const}$.   Therefore the tomography problem (\ref{eq:2.37})  is 
exactly the same as in the time-independent case.

In the case when $\gamma$  does not hit  obstacles,  (\ref{eq:2.37})
reduces to
\begin{equation}                            \label{eq:2.38}
\exp\{i\int_{-\infty}^\infty[A^{(3)}(x_0^{(0)}+s\omega,t_0)-
A^{(1)}(x_0^{(0)}+s\omega,t_0)]\cdot\omega ds\}=1
\end{equation}
We took in (\ref{eq:2.38})  $s_0=-\infty,\ s_1=+\infty$.
 
If conditions of Theorem \ref{theo:1.1}  are satisfied  then 
(\ref{eq:2.38})
implies  (see  [E2], \S  2)  that there exists  
$\varphi(x,t)\in C^\infty(\overline{D})$
such that
\begin{equation}                          \label{eq:2.39}
A^{(1)}(x,t)-A^{(3)}(x,t)=\frac{\partial \varphi}{\partial x},\ \ \
\varphi|_{\partial\Omega_0\times[0,T]}=0.
\end{equation}
If conditions of Theorem 1.2  are satisfied  then (\ref{eq:2.37})
implies  (see [E3],  pages 1512-1515)  that  there exists   $c(x,t) \in  G(D)$  
such that
\begin{equation}                                \label{eq:2.40}
A^{(3)}-A^{(1)}=ic^{-1}\frac{\partial c}{\partial x},\ \ 
c|_{\partial\Omega_0\times[0,T]}=1.
\end{equation}
Using (\ref{eq:2.39})  or (\ref{eq:2.40})
we make a gauge transformation $u_4=c^{-1}u_3$,  where $c=e^{i\varphi}$
in the case (\ref{eq:2.39}).  Then   $(i\frac{\partial}{\partial t}- H_4)u_4=0$,
where $i\frac{\partial}{\partial t}-H_4$  is the operator  with electromagnetic potentials 
$A^{(4)}=A^{(3)}-ic^{-1}\frac{\partial c}{\partial x},\ 
V^{(4)}=V^{(3)}+ic^{-1}\frac{\partial c}{\partial t}$.   Note that 
$A^{(1)}=A^{(4)}$.  Since  $c=1$  on  $\partial\Omega_0\times(0,T)$  we have
that  $\Lambda_4=\Lambda_3=\Lambda_1$,  where $\Lambda_4$  is the D-to-N 
operator corresponding  to  $i\frac{\partial}{\partial t}-H_4$. 
Apply the Green's formula (\ref{eq:2.35})  to 
$(i\frac{\partial}{\partial t}-H_1)u_1=0$  and
$(i\frac{\partial}{\partial t}-H_4^*)u_5=0$,  where  $u_5(x,T)=0,\ 
u_5|_{\partial\Omega_0\times(0,T)}=g_1,\ u_5|_{\partial\Omega'(t)}=0,\ 
\forall t\in[0,T]$.
 
Since  $A^{(1)}=A^{(4)}$ 
and $\Lambda_1=\Lambda_4$  we get
\begin{equation}                              \label{eq:2.41}
\int_D(V^{(1)}(x,t)-V^{(3)}(x,t)-ic^{-1}\frac{\partial c}{\partial t})
u_1\overline{u_5}dxdt=0.
\end{equation}
Substituting geometric optic solutions (\ref{eq:2.32}) 
and (\ref{eq:2.36})  with $H_3^*$ replaced by $H_4^*$  we get
\begin{equation}                              \label{eq:2.42}
\sum_{j=1}^r\int_{\gamma_j}(V^{(1)}-V^{(3)}-ic^{-1}\frac{\partial c}{\partial t})ds=0.
\end{equation}
Note that (\ref{eq:2.42})  holds for each  $t\in[0,T]$  and
for each broken ray in $D_t$.   It shown in [E2], \S 3,  that 
(\ref{eq:2.42})  implies that
\begin{equation}                             \label{eq:2.43}
V^{(1)}=V^{(3)}+ic^{-1}\frac{\partial c}{\partial t}
\end{equation}
Therefore (\ref{eq:2.40}),
(\ref{eq:2.43})
prove Theorem \ref{theo:1.2}.

Analogously,  when $\gamma$ does not hit the obstacles we have,
instead of (\ref{eq:2.42}):
\begin{equation}                             \label{eq:2.44}
\int_\gamma(V^{(1)}-V^{(3)}-ic^{-1}\frac{\partial c}{\partial t})ds=0.
\end{equation}
It is shown in [E2], \S2,  that (\ref{eq:2.44})
also implies (\ref{eq:2.43})  with  $c=e^{i\varphi(x,t)}$.  
Therefore Theorem \ref{theo:1.1} also holds.

{\bf Remark 2.2.}
Deriving formulas of the form (\ref{eq:2.37})  
we assumed  that  $\gamma=\gamma_1\cup...\cup\gamma_r$  are broken
nontangential rays that do not have focal  points  (i.e. do not intersect
the caustics set).  It was  noted in [E4],  page 29, 
 that these formulas still hold when $\gamma$ contains a generic caustics point.

Consider,  for simplicity,  the case $n=2$.  Let $(x_{10},x_{20},t_0)$  be the point
of the intersection of $\gamma$  with the caustics set.  In the
neighborhood of $(x_{10},x_{20},t_0)$  we replace the ansatz
$\sum_{p=0}^N\frac{a_p(x,t)}{(ik)^P}e^{-ik^2t+ik\psi(x,t)}$
with the following ansatz  (c.f. [V]):
\begin{equation}                      \label{eq:2.45}
v=\sum_{p=0}^N\int_{-\infty}^\infty\frac{b_p(x,t,\xi_1)}{(ik)^{p-\frac{1}{2}}}
e^{-ik^2t+i\varphi(x_2,t,\xi_1)-ix_1\xi_1}d\xi_1.
\end{equation}
At the caustics set we have:
\[
\varphi_{\xi_1}-x_1=0,\ \ \varphi_{\xi_1^2}=0
\]
We shall assume that $\varphi_{\xi_1^3}\neq 0$  on this set.
Then the following estimate holds in a neighborhood of $(x_{10},x_{20},t_0)$:
\begin{equation}                      \label{eq:2.46}
|v|\leq\frac{Ck^{\frac{1}{6}}}{1+k^{\frac{1}{6}}d^{\frac{1}{4}}}
\leq\frac{C}{d^{\frac{1}{4}}},
\end{equation}
where $d(x,t)$ is the distance to the caustics set.  It was shown in [E4]  
that the estimate (\ref{eq:2.46}) leads to the proof of 
(\ref{eq:2.37})  and(\ref{eq:2.42}).

\section{The Schr\"{o}dinger operator with time-dependent Yang-Mills 
potentials.} 
\label{section 3}
\init

The Schr\"{o}dinger equation with time-dependent Yang-Mills potentials 
has the following form:
\begin{equation}                                    \label{eq:3.1}
Lu\stackrel{def}{=}i\frac{\partial u(x,t)}{\partial t}-
\sum_{j=1}^n\left(-iI_m\frac{\partial}{\partial x_j}-A_j(x,t)\right)^2u(x,t)
+V(x,t)u(x,t),
\end{equation}
where  
$A_j,\ 1\leq j\leq n,\ V(x,t), u(x,t)$  are $m\times m$  matrices,
$I_m$  is the identity matrix in $\C^m$.  We assume that 
$A_j,\ 1\leq j\leq n, \ V$  are smooth
in $\overline{\Omega}_0\times[0,T]$  and without 
lost of generability (see Remark 3.1) we can assume that 
$A(x,t)=0, V(x,t)=0$ near $\partial\Omega_0\times(0,T)$.
We also  assume  that $A_j(x,t)$  are self-adjoint matrices,  $1\leq j\leq n$.
We consider (\ref{eq:3.1}) in $\Omega_0\times(0,T)$  with the initial
condition  
(\ref{eq:1.2})  and  the boundary condition $u|_{\partial\Omega_0\times(0,T)}=f$.
We do not consider domains with 
obstacles in the case of equation (\ref{eq:3.1}).

The gauge group $G$  now is the group of  $m\times m$ nonsingular smooth 
matrices          $g(x,t)$ in $\overline{\Omega}_0\times[0,T]$.
Let 
$L^{(p)}$ be two  Schr\"{o}dinger operators of the form (\ref{eq:3.1})
with Yang-Mills potentials $A^{(p)}(x,t),V^{(p)}(x,t),p=1,2.$

Let 
\begin{equation}                                \label{eq:3.2}
\Lambda_pf_p=(\frac{\partial}{\partial \nu }
-i A\cdot \nu)u_p|_{\partial\Omega_0\times(0,T)}, 
\ \ p=1,2,
\end{equation}
be corresponding D-to-N operators,   where  $L^{(p)}u_p=0$  in
$\Omega_0\times(0,T),\ u_p=0$  when $t=0,\ p=1,2,\ 
u_p|_{\partial\Omega_0\times(0,T)}=f_p$.  If  $u_2=g^{-1}u_1,\ \ \ g\in G,$ 
then
\begin{eqnarray}                          \label{eq:3.3}
A_j^{(2)}=g^{-1}A_j^{(1)}g+ig^{-1}\frac{\partial g}{\partial x_j},\ \ \
1\leq j\leq n,
\\
V^{(2)}=g^{-1}V^{(1)}g-ig^{-1}\frac{\partial g}{\partial t}.
\nonumber
\end{eqnarray}
Yang-Mills potentials $(A^{(2)},V^{(2)})$ and $(A^{(1)},V^{(1)})$ are called gauge
equivalent.   As in the case of Theorem \ref{theo:1.1} we have:
\begin{theorem}                         \label{theo:3.1}
If the D-to-N operators $\Lambda_1$  and $\Lambda_2$  are gauge equivalent 
on $\partial\Omega_0\times (0,T)$  (c.f.  (\ref{eq:1.7}) ) then 
$(A^{(1)},V^{(1)})$ and $(A^{(2)},V^{(2)})$  are gauge equivalent in
$\overline{\Omega}_0\times[0,T]$.
\end{theorem}

The start of the proof of Theorem \ref{theo:3.1} is the same as of 
Theorem \ref{theo:1.1}.  We construct  the geometric optics solution of 
(\ref{eq:3.1})  of the form:
\begin{equation}                       \label{eq:3.4}
u_N=e^{-ik^2t+ik(x,\omega)}\sum_{p=0}^N\frac{1}{(ik)^p}a_p(x,t,\omega)+u^{(N+1)},
\end{equation}
where $a_0(x,t)=\chi_1(t)\chi_2(\tau)c(x,t,\omega)$  and
\begin{eqnarray}                        \label{eq:3.5}
\omega\cdot\frac{\partial c}{\partial x}  - iA(x,t)\cdot\omega c(x,t)=0,\ \ 
s> s_0,
\\
c|_{s=s_0}=I_m,
\nonumber
\end{eqnarray}
$a_p(x,t,\omega)$  are solutions of the equations
$2\omega\cdot(-i\frac{\partial}{\partial x}-A)a_p=
La_{p-1}$  for  $s>s_0,\ a_p|_{s=s_0}=0,\ p\geq 1$.

Here $s,s_0,\chi_1(t),\chi_2(\tau)$  are the same as in
(\ref{eq:2.4}),  (\ref{eq:2.5}).  Note  that
since $A(x,t)$ have a compact support  we can take $s_0=-\infty$.
Finally, $u^{(N+1)}$
satisfies
\begin{eqnarray}                       \label{eq:3.6}
Lu^{(N+1)}=\frac
{-e^{-ik^2t+ik(\omega\cdot x)}}
{(ik)^N}L a_N,\ \ (x,t)\in \Omega_0\times (0,T),
\\
u^{(N+1)}(x,0)=0,\ \ u^{(N+1)}|_{\partial\Omega_0}\times (0,T)=0.
\nonumber
\end{eqnarray}
The existence of $u^{(N+1)}$  follows  from  Lemma
\ref{lma:2.1} that holds without changes in the case of equation
(\ref{eq:3.1}).   Since  $\Lambda_2=g_0^{-1}\Lambda_1g_0$  for some
$g\in G(\Omega_0\times(0,T))$  where $g_0=g|_{\partial\Omega_0\times(0,T)}$
we have that $\Lambda_3=\Lambda_1$   where $\Lambda_3$  is the D-to-N  operator 
corresponding 
to $L_3$  with potentials $A^{(3)}=g^{-1}A^{(2)}g-
ig^{-1}\frac{\partial g}{\partial x},\ 
V^{(3)}=g^{-1}V^{(2)}g+
ig^{-1}\frac{\partial g}{\partial t}$.
Substituting the geometric optics solutions for $L_1$ and $L_3^*$ into
the Greens formula of the form (\ref{eq:2.35})
we get as in \S2 (c.f.  (\ref{eq:2.38}) 
and also [E4],  pages 29-30):
\begin{equation}                             \label{eq:3.7}
c_{30}^{-1}(+\infty,y_2,t,\omega)c_{10}(+\infty,y_2,t,\omega)=I_m,
\end{equation}
where $y_1=x\cdot\omega,\ y_2=x-(x\cdot\omega)\omega,\ c_j(x,t,\omega),\ j=1,3,$
  are the solution  of the differential equation  of the form (c.f. (\ref{eq:3.5}) ):
\begin{equation}                             \label{eq:3.8}
\omega\cdot\frac{\partial c_j(x,t,\omega)}{\partial x} 
- i(A_j(x,t)\cdot\omega)c_j(x,t,\omega)=0,
\end{equation}
$c_{j0}(y_1,y_2,t,\omega)$  is $c_j(x,t,\omega)$  in $(y_1,y_2)$ coordinates,
$c_{j0}(-\infty,y_2,t,\omega)=I_m$.  The matrix $c_{j0}(+\infty,y_2,t,\omega)$  is
called the non-abelian Radon transform of $A^{(j)}(x,t),\ j=1,3$.

The following lemma was proven in [E5]  (see also [No]):
\begin{lemma}                      \label{lma:3.1}
 If the non-abelian Radon transform of $A^{(1))}$ and $A^{(3)}$  are equal
then $A^{(1))}$ and $A^{(3)}$ are gauge equivalent with 
$g\in G_0(\Omega_0\times[0,T]$,  i.e. 
$g=I_m$  on $\partial\Omega_0\times[0,T]$.
\end{lemma}
The proof of Lemma \ref{lma:3.1} is much harder than the similar result
in the abelian case  (i.e. $m=1$), and involves  a lot of the complex
analysis.  We consider briefly the most important case of two dimensions.  
The case $n\geq 3$ is reduced easily  to the case  $n=2$.   The main
ingredient of the proof is the study  of an equation  of the form (\ref{eq:3.8})
with complex parameters:
\begin{equation}                      \label{eq:3.9}
\zeta_1\frac{\partial c}{\partial x_1}+\zeta_2\frac{\partial c}{\partial x_2}
=i(A_1(x,t)\zeta_1+A_2(x,t)\zeta_2)c,
\end{equation}
where
$\zeta_1^2+\zeta_2^2=1,\ \zeta_i\in\C,\ i=1,2.$  Define
$\zeta_1(\tau)=\frac{1}{2}(\tau+\frac{1}{\tau}),\ 
\zeta_2(\tau)=\frac{i}{2}(\tau-\frac{1}{\tau})$,  
where $\tau\in\C\setminus \{0\}$.
It was proven in [ER2]  that there exists solutions $c_\pm(x,t,\tau)$ of
(\ref{eq:3.9}),  smooth in $(x,t,\tau)$  and such that $c_+(x,t,\tau)$
is analytic and a nonsingular matrix in $\tau$  when $|\tau|<1,\ c_-(x,t,\tau)$  
is analytic and nonsingular in $\tau$ for $|\tau|>1$.  Let $\tau=|\tau|e^{i\varphi}$.
When $|\tau|\rightarrow 1$  we have that $(\zeta_1(\tau),\zeta_2(\tau)\rightarrow
\omega(\varphi)\stackrel{def}{=}(\cos\varphi,-\sin\varphi)$.  Therefore
$c_\pm(x,t,e^{i\varphi})$  are solutions of equation 
of the form (\ref{eq:3.8})  with $\omega(\varphi)=(\cos\varphi,-\sin\varphi)$.
As before,  denote $y_1=(x\cdot\omega(\varphi))$  and 
$y_2=(x\cdot\omega_\perp(\varphi))$
where $\omega_\perp(\varphi)=(\sin\varphi,\cos\varphi)$.  Let 
$c_\pm^{(0)}(y_1,y_2,t,e^{i\varphi})$  be $c_\pm(x,t,e^{i\varphi})$  in
$(y_1,y_2)$  coordinates and let $c_0(y_1,y_2,t,\omega(\varphi))$ 
be the solution of the Cauchy problem (\ref{eq:3.8}).   We have:
\begin{equation}                     \label{eq:3.10}
c_0(y_1,y_2,t,\omega(\varphi))=
c_\pm^{(0)}(y_1,y_2,t,e^{i\varphi})
(c_\pm^{(0)}(-\infty, y_2,t,e^{i\varphi}))^{-1}
\end{equation}                       
by the uniqueness of the Cauchy problem for (\ref{eq:3.8}).  Putting
$y_1=+\infty$  we get
\begin{equation}                     \label{eq:3.11}
c_0(+\infty,y_2,t,\omega(\varphi))=
c_\pm^{(0)}(+\infty,y_2,t,e^{i\varphi})(c_\pm^{(0)}(-\infty, y_2,t,e^{i\varphi}))^{-1}.
\end{equation}  
It follows from (\ref{eq:3.7})  and  (\ref{eq:3.11})  that
\begin{eqnarray}                      \label{eq:3.12} 
c_{1,\pm}^{(0)}(+\infty,y_2,t,e^{i\varphi})
(c_{1,\pm}^{(0)}(-\infty,y_2,t,e^{i\varphi}))^{-1}
\\
=c_{3,\pm}^{(0)}(+\infty,y_2,t,e^{i\varphi})
(c_{3,\pm}^{(0)}(-\infty,y_2,t,e^{i\varphi}))^{-1},
\nonumber
\end{eqnarray}
where
$c_{j,\pm}(x,t,e^{i\varphi})$  correspond to the potentials $A^{(j)}(x,t),j=1,3$.
The continuation of the proof of Lemma \ref{lma:3.1}  when $n=2$ 
is the same as in [E5],
\S 2.

Now we shall consider the case $n\geq 3$ assuming that (\ref{eq:3.7}) holds. 
Let $\omega=(\omega_1,...,\omega_n)\in S^{n-1}$  and let $\omega_{12}=
(\omega_1,\omega_2,0,...,0)\in S^{n-1}$.  Consider equations  
(\ref{eq:3.8})
with $\omega$  replaced  by $\omega_{12}$.   It is a two-dimensional problem
with $x_3,...,x_n$  and  $t$ as parameters.   It follows  from Lemma \ref{lma:3.1}
for $n=2$ that there exists $g(x,t),\ \det g(x,t)\neq 0,\ g=I_m$ on
$\partial\Omega_0\times[0,T]$  such that
\begin{equation}                               \label{eq:3.13}
A_1^{(3)}=g^{-1}A_1^{(1)}g+ig^{-1}\frac{\partial g}{\partial x_1},
\  \ \ \ 
A_2^{(3)}=g^{-1}A_2^{(1)}g+ig^{-1}\frac{\partial g}{\partial x_2}
\end{equation}
for all $(x,t)$.   Replacing
$\omega_2$  by $\omega_j,\ j>2, $
we get that there exists $\tilde{g}\in G_0(\Omega_0\times[0,T])$  such
that 
\begin{equation}                               \label{eq:3.14}
A_1^{(3)}=\tilde{g}^{-1}A_1^{(1)}\tilde{g}
+i\tilde{g}^{-1}\frac{\partial \tilde{g}}{\partial x_1},
\  \ \ \ 
A_j^{(3)}=\tilde{g}^{-1}A_j^{(1)}\tilde{g}
+i\tilde{g}^{-1}\frac{\partial \tilde{g}}{\partial x_2}
\end{equation}
for all $(x,t)$.
Comparing (\ref{eq:3.13})  and
(\ref{eq:3.14})  we have that $g$ and $\tilde{g}$  satisfy the same system 
of differential equations 
\[
i\frac{\partial g}{\partial x_1}+A_1^{(1)}g - gA_1^{(3)}=0,\ \ \ 
g|_{\partial\Omega_0\times[0,T]}=I_m
\]
for each 
 $(x_2,...,x_n,t)$.

By the uniqueness theorem for the differential equations we get that 
$\tilde{g}=g$.  Therefore
\begin{equation}                               \label{eq:3.15}
A_j^{(3)}=g^{-1}A_j^{(1)}g+ig^{-1}\frac{\partial g}{\partial x_j},
\  \ \ \ j\geq 3.
\end{equation}
Denote
\begin{equation}                                \label{eq:3.16}
A^{(4)}=gA^{(3)}g^{-1}-i\frac{\partial g}{\partial x}g^{-1},\ \ \ 
V^{(4)}=gV^{(3)}g^{-1}+i\frac{\partial g}{\partial x}g^{-1}.
\end{equation}
Then $A^{(4)}(x,t)=A^{(1)}(x,t)$.

Now we shall show that $V^{(4)}=V^{(1)}$,  i.e.  $(A^{(3)},V^{(3)})$  and
$(A^{(1)},V^{(1)})$
are gauge equivalent.  As before,  it is enough to consider the case
$n=2$.  Using the Green's formula for $L^{(1)}$ and $L^{(4)*} $
we get,  as in \S 2 (c.f. (\ref{eq:2.41}),  (\ref{eq:2.44}) ), that  
\begin{equation}                             \label{eq:3.17}
\int_{-\infty}^\infty  c_{10}^{-1}(y_1,y_2,t,\omega)(V^{(1)}(x,t)
-V^{(4)}(x,t))c_{10}(y_1,y_2,t,\omega)dy_1=0.
\end{equation}
Here $y_1=x\cdot\omega,\ y_2=x-(x\cdot\omega)\omega$.
Substituting (\ref{eq:3.10}) into (\ref{eq:3.17}],  multiplying (\ref{eq:3.16})
by $c_\pm^{-1}(-\infty,y_2,t,e^{i\varphi})$ from the left and
$c_\pm (-\infty,y_2,t,e^{i\varphi})$ from the right,
we get
\begin{equation}                             \label{eq:3.18}
\int_{-\infty}^\infty  c_\pm^{-1}(y_1,y_2,t,e^{i\varphi})(V^{(1)}-V^{(4)})
c_\pm(y_1,y_2,t,e^{i\varphi})dy_1=0.
\end{equation}

The continuation of proof that (\ref{eq:3.18})
leads to $V^{(1)}-V^{(4)}=0$  is the same as in [E1],  pages 61-62.
\qed

{\bf  Remark 3.1.}  In this remark we shall show that without loss
of generality one can assume that $A(x,t)$  and $V(x,t)$ are smooth and 
equal to zero near $\partial\Omega_0\times[0,T] $  (c.f.  [NSU], [ER1], [E1]).
Let $L_p$ be two operators in $\Omega_0\times[0,T]$  such that  
$\Lambda_1=\Lambda_2$  on   $\partial\Omega_0\times(0,T)$.  As in [E3],
pages 53-54,  one can replace $(A^{(2)},V^{(2)})$  in $\Omega_0\times[0,T]$
by gauge equivalent $(A^{(3)}, V^{(3)})$  with the gauge $g\in G_0$  such 
that $A^{(3)}(x,t)\cdot n(x,t)=A^{(1)}(x,t)\cdot n(x,t)$  
near $\partial\Omega_0\times[0,T]$,  where $n(x,t)$  is
a smooth extension near $\partial\Omega_0\times[0,T]$  of the unit normal 
vector $n(x)$ on $\partial\Omega_0$.
We shall show that $\Lambda_1=\Lambda_3$  on $\partial\Omega_0\times(0,T)$
implies that $A^{(1)}=A^{(3)}$  and $V^{(1)}=V^{(3)}$ on 
$\partial\Omega_0\times[0,T]$  with its derivatives.   
Introduce in a neighborhood $U$  of 
an arbitrary point $(x_0,t_0)\in\partial\Omega_0\times[0,T_0]$
a system of coordinates $(y,t)$  such that the equation of 
$\partial\Omega_0$  is  $y_n=0$.   Let
$(\eta_1,...,\eta_n,\eta_0)$  be the dual coordinates to $(y,t)$
in the cotangent space $T^*(U\times[0,T))$.  The Schr\"{o}dinger 
operators $L_1$ and $L_3$ are elliptic microlocally in the region
where $|\eta_0|<\e\sum_{k=1}^n\eta_k^2$  (c.f. [E6],  Remark 2.2).
Therefore the elliptic parametrix  for $\Lambda_1$  and $\Lambda_3$  applies 
in this region.   Therefore,  as in [E1],  Proposition 2.1,  we get
that
\[
\frac{\partial^{p+q}A^{(1)}}{\partial x^p\partial t^q}= 
\frac{\partial^{p+q}A^{(3)}}{\partial x^p\partial t^q}
\ \ \ \mbox{and}\ \ \ 
\frac{\partial^{p+q}V^{(1)}}{\partial x^p\partial t^q}=
\frac{\partial^{p+q}V^{(3)}}{\partial x^p\partial t^q}
\]
on $\partial\Omega_0\times(0,T)$  for any  $p$ and $q$.
Therefore we can extend smoothly  $(A^{(1)},V^{(1)})$  and
$(A^{(3)},V^{(3)})$ 
from $\overline{\Omega}_0\times[0,T]$  to a large domain
$\overline{\tilde{\Omega}}_0\times[0,T]$  in a such way that
$A^{(1)}=A^{(3)},\ V^{(1)}=V^{(3)}$  in
($\overline{\tilde{\Omega}}_0\setminus\Omega_0)\times[0,T]$ 
and $(A^{(1)},V^{(1)})=
(A^{(3)},V^{(3)})=(0,0)$ 
near $\partial\tilde{\Omega}_0\times[0,T]$ .
It remains to show that $\Lambda_1=\Lambda_3$  on 
$\partial\Omega_0\times(0,T)$  implies
that $\Lambda_1=\Lambda_3$  on $\partial\tilde{\Omega}_0\times[0,T]$.
Take any smooth $f$ on $\partial\tilde{\Omega}_0\times[0,T]$  and consider 
the initial-boundary value problem
\begin{eqnarray}                   
\nonumber
L^{(1)}u_1=0\ \ \ \ \ \mbox{in}\ \ \tilde{\Omega}_0\times(0,T),
\\
\nonumber
u_1(x,0)=0,\ \ \ \ \ \ \ u_1|_{\partial\tilde{\Omega}_0\times[0,T]}=f
\end{eqnarray}
Let $f_1$ be the restriction of $u_1$ to $\partial\Omega_0\times(0,T).$
Consider the initial-boundary value problem
\begin{eqnarray}                   
\nonumber
L_3u_3=0\ \ \ \ \ \mbox{in}\ \ \Omega_0\times(0,T),
\\
\nonumber
u_3(x,0)=0,\ \ \ \ \ \ \ u_3|_{\partial\Omega_0\times[0,T]}=f_1.
\end{eqnarray}
Let $w_3=u_1$  in 
$(\overline{\tilde{\Omega}}_0\setminus\Omega_0)\times(0,T),\ w_3=u_3$
in $\Omega_0\times(0,T)$.  Since $\Lambda_1=\Lambda_3$  on
$\partial\Omega_0\times(0,T)$ we get that $w_3$   is the solution of $L_3w=0$
in $\tilde{\Omega}_0\times(0,T)$.  
Moreover,  $w_3(x,0)=0,\ x\in\tilde{\Omega}_0$
and $w_3|_{\partial\tilde{\Omega}\times(0,T)}=f$.  Since  $L^{(1)}=L^{(3)}$
and  $u_1=w_3$  in $(\tilde{\Omega}_0\setminus\Omega_0)\times(0,T)$
we have that $\Lambda_1=\Lambda_3$  on $\partial\tilde{\Omega}_0\times(0,T)$.
Therefore we have reduced the inversee problem in $\Omega_0\times(0,T)$ to
an equivalent inverse problem in a larger domain $\tilde{\Omega}_0\times(0,T)$.
Note that the same proof works when $\Omega_0\times(0,T)$ contains obstacles.

\section{Examples of the gauge equivalence classses and the Aharonov-Bohm effect.}
\label{section 4}
\init

The results of  \S 2 are valid in the multiconnected domains $D$.  
In such domains there are electromagnetic potentials $(A(x,t),V(x,t))$  and
$(A'(x,t),V'(x,t))$ that correspond to the same electromagnetic fields in
$D : \mbox{curl\ } A=\mbox{curl\ } A'=B,\ E=-\frac{\partial V}{\partial x}-
\frac{\partial A}{\partial t}=
-\frac{\partial V'}{\partial x}-\frac{\partial A'}{\partial t},$
however $(A,V)$  and $(A',V')$  are not gauge equivalent.
There is a simple description of all gauge equivalent classes of
potentials in $D$ (c.f. [WY], [OP], [Va]):

Fix a point $(x_0,t_0)\in\partial\Omega_0\times[0,T]$  and denote
by $\mathcal{P}$  the group
of all paths $\gamma$  in $\overline{D}$  that start and end  at $(x_0,t_0)$.  
Let $\int_\gamma  A(x,t)\cdot dx-V(x,t)dt$  be a line integral in $D$.  Denote
\begin{equation}                               \label{eq:4.1}
R(A,V,\gamma)=\exp\{i(\int_\gamma -A(x,t)\cdot dx+V(x,t)dt)\}.
\end{equation}
$R$ is called the nonintegrable phase factor (c.f. [WY],  [OP]).
$R(A,V,\gamma_0)$ is well-defined for any closed path $\gamma_0$ in
$\overline{D}$ :
Let $\gamma_1$  be a path
connecting $(x_0,t_0)$  with an arbitrary point $(x_1,t_1)\in \gamma_0$
and let $\gamma\in\mathcal{P}$  be a close path  $\gamma_1\cup\gamma_0$.  
It is obvious that $R(A,V,\gamma)=R(A,V,\gamma_0)$.
If $(A,V)$  and $(A',V')$  are gauge equivalent  then
\begin{eqnarray}                      
\nonumber
R'R^{-1}=\exp\left\{i \left(\int_\gamma(A-A')\cdot dx + (V'-V)dt\right)\right\}
\\
=\exp\left\{-i\left(\int_\gamma ic^{-1}\frac{\partial c}{\partial x}\cdot
dx + ic^{-1}\frac{\partial c}{\partial t}dt\right)\right\}=
\exp\left\{\int_\gamma d (\log c)\right\}=1,
\nonumber
\end{eqnarray}
since  
$c(x,t)\in G(D)$  is a single-valued in $\overline{D}$.  Therefore
$R'=R$.  And vice versa,   if $R(A,V,\gamma)=R(A',V',\gamma)$  for 
all $\gamma\in \mathcal{P}$  then
$(A,V)$  and $(A',V')$  are gauge equivalent :  Let
\begin{equation}                            \label{eq:4.2}
\exp\left\{i \left(\int_\gamma(A'-A)\cdot dx + (V-V')dt\right)\right\}=1
\end{equation}
for all $\gamma\in \mathcal{P}$.  
Let $(x,t)$  be an arbitrary point in $\overline{D}$.
Denote by $c(x,t)$  the integral 
\begin{equation}                            \label{eq:4.3}
c(x,t)=
\exp\left\{i \left(\int_{\gamma(x,t)}(A'-A)\cdot dx + (V-V')dt\right)\right\}
\end{equation}
where $\gamma(x,t)$  is an arbitrary path in $\overline{D}$ connecting
$(x_0,t_0)$ and $(x,t)$.  It follows from (\ref{eq:4.2})   that
(\ref{eq:4.3}) does not depend on $\gamma(x,t)$  connecting 
$(x_0,t_0)$  and $(x,t)$.   Differentiating (\ref{eq:4.3})  in $x$
and $t$  we get that $(A,V)$  and $(A',V')$  are gauge equivalent.  

For fixed $(A,V)\ \ R(A,V,\gamma)$  defines a map of the group of
paths $\mathcal{P}$  to the group  $U$  of complex numbers $z$
with absolute value 1,  i.e.  $z=e^{i\alpha},\ \alpha\in\R $  (we assume 
in this section that $(A,V)$ are real-valued).  The image of this map
$U(A,V)$  is a subgroup of $U$ and it is called the monodromy group
of the connection $\mathcal{A} =\sum_{j=1}^nA_jdx_j-Vdt$  (c.f. [Va]).
We will often write $\mathcal{A}$  instead of $(A,V)$.
Therefore there is one-to-one correspondence between the gauge 
equivalence classes 
of electromagnetic potentials and the monodomy groups.  It was shown 
by Aharonov and Bohm (c.f. [AB])  that potentials belonging to
different gauge equivalence classes  produce a different physical impact
and it can be detected in experiments.  This phenomenon is called 
the Aharonov-Bohm effect.

Let $S$  be a two-dimensional smooth surface  in $D$  such that 
$\partial S=\gamma$.  It is convenient to represent the 
electromagnetic field $(E,B)$  as a two-form
$\mathcal{F}=\sum_{1\leq j<k\leq n}B_{jk}dx_j\wedge dx_k+
\sum_{k=1}^n E_k dx_k\wedge dt$,  where $B_{jk}=-B_{kj}$.
Note that $\mathcal{F}=d\mathcal{A}$.
Therefore $B_{kj}=\frac{\partial A_j}{\partial x_k}-
\frac{\partial A_k}{\partial x_j},
\ E_j=-\frac{\partial A_j}{\partial t}-\frac{\partial V}{\partial x_j},\ 
1\leq j\leq n,\ 1\leq k\leq n$.
When $n=3$  we have $B=(B_1,B_2,B_3),\ E=(E_1,E_2,E_3),\ 
B_1=B_{23},\ B_2=-B_{13}, \ B_3=B_{12}$.   Therefore,  when $n=3$  
\begin{equation}                            \label{eq:4.4}
B=\mbox{curl\ }A,\ \ \ \ E=-\frac{\partial A}{\partial t}
-\frac{\partial V}{\partial x}
\end{equation}
Using the Stoke's formula $\int_\gamma \mathcal{A}=
\int_S d\mathcal{A}$ we get (c.f. [OP]):
\begin{eqnarray}                        \label{eq:4.5}
\int_\gamma A\cdot dx - Vdt=\int_S(B_1dx_2dx_3+B_2dx_1dx_3+B_3dx_1dx_2
\\
+E_1dtdx_1+E_2dtdx_2+E_3dtdx_3).
\nonumber
\end{eqnarray}
The right hand side of (\ref{eq:4.5})  is called the electromagnetic 
flux.  Note that the electromagnetic flux does not depend on the surface
$S$ such that $\partial S=\gamma$.
This follows from the three-dimensional Stoke's theorem 
and the equality $d\mathcal{F}=d(d\mathcal{A})=0$.

In the case $n=2$  formula (\ref{eq:4.5}) for the electromagnetic flux
has the form:
\begin{equation}                         \label{eq:4.6}
\int_\gamma A_1dx_1+A_2dx_2-Vdt=
\int_S(B_3dx_1dx_2+E_1dtdx_1+E_2dtdx_2).
\end{equation}
Formulas (\ref{eq:4.5}),  (\ref{eq:4.6})  can be used for the computations 
of $R(A,V,\gamma)$,
especially when either the magnetic or electric fields is "shielded"
inside the obstacles,  i.e.  when either $B$ or $E$ is nonzero in
$\Omega'(t)=\cup_{j=1}^m\Omega_j(t)$ and zero in $D$.

The following examples give a description of gauge equivalence classes of
electromagnetic potentials using the electromagnetic fluxes.

Consider first an example when electromagnetic field is time-independent. 
Let $n=3,\ D_0=\{x:|x|<r\}$  and the single obstacle $D_1$ is a torus
inside $D_0$  as in Remark 3.3 in [E3].  Let $B(x)$  be a magnetic
field in $D_0$  such that $\mbox{supp\ }B\subset D_1$,  
 i.e. $B$ is "shielded"  in  $D_1$.
Denote by $b_0$  the flux of $B$ over arbitrary cross-section $S_0$
of the torus $D_1$.

Note that $b_0$  does not depend on the choice of the
 cross-section since $\mbox{div\ }B=0$.
Let $A$ be the magnetic potential in $D_0$,  i.e. $\mbox{curl\ }A=B$
in $D_0$.  Since $B=0$  in $D_0\setminus D_1$  the integral  
$\int_\gamma A\cdot dx$
depends only  on the homotopy class of $\gamma$ in $D_0\setminus D_1$.
If $B'$  is another magnetic field in $D_0$  shielded in $D_1$ and 
$\mbox{curl\ }A'=B'$  in $D_0$  then potentials $(A,V)$  and
$(A',V')$ are gauge equivalent in $D_0\setminus D_1$ 
iff $V'=V$  and there exists $m_0\in \Z$  such that 
the flux $\int_{S_0}(B'\cdot n)ds = b_0+2\pi m_0$,  where $n$  is 
the unit vector normal to $S_0$.

In the case of time-dependent electromagnetic potential the situation is more
complicated 
because both the magnetic and the electric field can not be shielded inside
the obstacle.   Moreover,  we have to use the flux (\ref{eq:4.5}) 
instead of the magnetic flux only.

Consider  the case  $n=2,\ \Omega_0=\{x:|x|<r\}$  with an obstacle 
$\Omega_1$  moving with speed $v_0$
in the direction of $x_1$-axis,  $\Omega_1=\cup_{t\in(0,T)}\Omega_1(t)$,
where $\Omega_1(t)=\{(x,t):(x_1-v_0t)^2+x_2^2\leq r_1\},\ t\in[0,T]\},
\ r_1<r $  is small.
Here $D=(\Omega_0\times(0,T))\setminus\Omega_1$.
In this case $G(D)$  consists of the functions $c(x,t)$ having the form
\begin{equation}                           \label{eq:4.7}
c(x,t)=e^{im\theta+i\psi(x,t)},
\end{equation}
where $m\in \Z,\ \psi(x,t)\in C^\infty(\overline{D})$  and $\theta$ is
the angle in polar coordinates centered at $(v_0t,0)$  with $x_2=0,\ 
x_1-v_0t>0$  as the polar axis.  Consider an electric field $E$  "shielded"
by the "moving"  obstacle $\Omega_1$.
Suppose $E=(0,E_2)$  where
\begin{equation}                        \label{eq:4.8}
E_2(x,t)=\delta(x_1-v_0t)\delta(x_2)e(t),
\end{equation}
$e(t)\in C^\infty, \ e(t)=\mbox{const}$  when $t\not\in [t_1,t_1+\e]$,  where
$t_1>0,\ t_1+\e<T$.
Let $\Omega_1(x_{10})$  be the intersection of the plane $x_1=x_{10}$
with $\Omega_1$  and let $\gamma(x_{10})$  be a simple close curve
in $(\Omega_0\times(0,T))\cap \{x_1=x_{10}\}$  encircling $\Omega_1(x_{10})$.
By formula (\ref{eq:4.6}) we have 
\begin{equation}                       \label{eq:4.9}
\int_{\gamma(x_{10})}A_2 dx_2-     Vdt
=\int_{\Omega_1(x_{10})}E_2dx_2dt=\frac{1}{v_0}e(\frac{x_{10}}{v_0}).
\end{equation}
The right hand side of (\ref{eq:4.9})  is the electric flux over
$\Omega_1(x_{10})$  and $V(x,t)$ is an electric potential corresponding
to $E=(0,E_2)$.  Note  that  (\ref{eq:4.9})
holds with $e(t)$  replaced by 0 for any closed Jordan curve
in $\{x_1=x_{10}\}\setminus\Omega_1$.
Let $B_3(x,t)$  be the magnetic field caused by the electric  field $(0,E_2)$.
It follows from the Maxwell's equations that 
\[
\frac{\partial E_2}{\partial x_1}=\frac{\partial B_3}{\partial t}.
\]
We get from (\ref{eq:4.8})
\begin{equation}                      \label{eq:4.10}
B_3(x,t)=-\frac{1}{v_0}\delta(x_1-v_0t)\delta(x_2)e(t)-\frac{1}{v_0^2}
\theta(x_1-v_0t)\delta(x_2)\frac{\partial e(\frac{x_1}{v_0})}{\partial t},
\end{equation}
where $\theta(\tau)=1$  when $\tau>0,\ \theta(\tau)=0$  when $\tau<0$.

Let $\gamma(t)$  be a closed simple curve in 
$(\Omega_0\cap\{t_1=\mbox{const}\}) \setminus\Omega_1(t)$.  We have by 
(\ref{eq:4.6})
\[
\int_{\gamma(t)} A\cdot dx =\int_{S(t)}B_3dx_1dx_2 =-\frac{1}{v_0}e(t)
-\frac{1}{v_0^2}\int_{S(t)}\theta_1(x_1-v_0t)\delta(x_2)
\frac{\partial e(\frac{x_1}{v_0})}{\partial t} dx_1dx_2,
\]
where  $\partial S(t)=\gamma (t)$.  Since 
$\frac{\partial}{\partial t}e(\frac{x_1}{v_0})
=v_0\frac{\partial}{\partial x_1} e(\frac{x_1}{v_0})$  we get
\[
\int_{S(t)}\theta(x_1-v_0t)\delta(x_2)
\frac{\partial e(\frac{x_1}{v_0})}{\partial t}dx_1d_2
=v_0\left(e(\frac{\tilde{x}_1(t)}{v_0})-e(t)\right),
\]
where $\tilde{x}_1(t)$ is the intersection of the line 
$x_1=v_0t,\ x_2=0$ with $\gamma(t)$.
Therefore 
\begin{equation}                         \label{eq:4.11}
\int_{S(t)}B_3dx_1dx_2 =-\frac{1}{v_0}e(\frac{\tilde{x}_1(t)}{v_0}).
\end{equation}
Note that knowing the electric flux (\ref{eq:4.9}) for all $x_{10}$
we can determine the magnetic flux (\ref{eq:4.11}).
Let $E_2'(x,t)=\delta(x_1-v_0t)\delta(x_2)e'(t)$ be 
shielded by  the obstacle
$\Omega_1$, where $e'(t)=\mbox{const}$
when $t\notin[t_1,t_1+\e]$,
 and let $B_3'$  be the coresponding magnetic field.  Let $(A,V)$
and $(A',V')$  be any corresponding electromagnetic potentials  in
$\Omega_0\times(0,T)$.  If $(A,V)$
and $(A',V')$ are gauge equivalent in $D$ then there exists $m\in\Z$ such that
\begin{equation}                        \label{eq:4.12}
\frac{1}{v_0}e'(\frac{x_1}{v_0})=
\int_{\Omega_1(x_1)}E_2'dx_2dt=\int_{\Omega_1(x_1)}E_2dx_2dt+2\pi m=
\frac{1}{v_0}e(\frac{x_1}{v_0}) +2\pi m
\end{equation}
for all $x_1$.  Vice versa,  suppose (\ref{eq:4.12})  holds.  Then for the 
magnetic fluxes we also have
\[
\int_{S(t)}B_3'dx_1dx_2=\int_{S(t)}B_3dx_1dx_2+2\pi m
\]
for any $S(t)$ such that  $\gamma(t)=\partial S(t)$ is 
simple closed curve in $\Omega_0\setminus \Omega_1(t)$.
Moreover,  (\ref{eq:4.10}) and (\ref{eq:4.12}) imply that
$B_3-B_3'=0$ in $\Omega_0\setminus\Omega_1(t)$ for all $t$.
This means that $B_3-B_3'$ as $E_2-E_2'$ are confined to the obstacle $\Omega_1$.
Therefore for any close path $\gamma$ in $(\Omega_0\times(0,T))\setminus\Omega_1$
the line integral $\int_\gamma(A'-A)\cdot dx+(V-V')dt$ depends only on the
homotopy class of $\gamma$.  Making a homotopy of $\gamma$ to a multiple of
the $\gamma(x_{10})$ we get,  using (\ref{eq:4.12}),  that 
$\int_\gamma(A'-A)\cdot dx+(V-V')dt=2\pi m'$,
where $m'\in \Z$.
Therefore $(A,V)$
and $(A',V')$ are gauge equialent.

Finally, let $n=2$ and $\Omega_0,\ \Omega_1$ be  the same as above.
Consider an example of the electromagnetic fluxes caused by the magnetic field
$B_3(x,t)=\delta(x_1-v_0t)\delta(x_2)b(t)$ in $\Omega_0\times(0,T)$  
shielded in $\Omega_1$.
In this case we can allow $v_0=0$, i.e. the obstacle $\Omega_1$ may not move.
Let $\gamma(t_0)$  be a closed simple curve in $\Omega_0\times\{t=t_0\}$ that 
encircles $\Omega_1(t_0)$.  Then by (\ref{eq:4.6})
\[
\int_{\gamma(t_0)} A\cdot dx =\int_{\Omega_1(t_0)}B_3dx_1dx_2=b(t_0).
\]
Here $A(x,t)$  are arbitrary magnetic potentials in $\Omega_0\times(0,T)$  such
that $\mbox{curl\ }A=B(x,t)$.  For example,  we may take
\[
A=(A_1(x),A_2(x))=\frac{b(t)}{2\pi}\left(\frac{-x_2}{(x_1-v_0t)^2+x_2^2},
\frac{x_1-v_0t}{(x_1-v_0t)^2+x_2^2}\right).
\]
Let $B_3'(x,t)$  be an arbitrary magnetic field shielded in $\Omega_1$
and such that
\begin{equation}                     \label{eq:4.13}
\int_{\Omega_1(t)}B_3'(x,t)dx_1dx_2=b(t)+2\pi m
\end{equation}
for some $m\in\Z$.  Denote by $A'$ an arbitrary
magnetic potential in $\Omega_0\times(0,T)$  such that $\mbox{curl\ }A'=B_3'$
in $\Omega_0\times(0,T)$.
Then we have
\[
\exp(i\int_{\gamma(t)}A\cdot dx)=\exp(i\int_{\gamma(t_0)}A'\cdot dx)
\]
for any closed path in $(\Omega_0\times\{t=t_0\})\setminus\Omega_1(t_0),\ 
\forall t_0\in[0,T]$.
Therefore there exists $c(x,t)\in G(D)$  such that 
$A'-A=-ic^{-1}\frac{\partial c}{\partial x}$  in $D$.

Define $V'=V+ic^{-1}\frac{\partial c}{\partial t}$.  Then 
$(A',V')$  is gauge equivalent
to $(A,V)$.  Vice versa,  if $(A',V')$  and $(A,V)$  are gauge 
equivalent then (\ref{eq:4.13})  holds.

Note that in this example the magnetic field  is shielded in
$\Omega_1$  but the electric field $E=-\frac{\partial V}{\partial x}
-\frac{\partial A}{\partial t}$  is not shielded in $\Omega_1$.
In the previous example the electric field $E$  was shielded in $\Omega_1$
but the magnetic field was not.

\end{document}